\documentclass[11pt,twoside]{amsart}

\usepackage{amsfonts}

\title[Triangulated version of Green correspondence for adjoint functors]{Remarks on a triangulated version of Auslander-Kleiner's Green correspondence}

\usepackage[all]{xy}
\usepackage{graphicx}
\usepackage{yhmath}
\usepackage{mathdots}
\usepackage{MnSymbol}

\author{Alexander Zimmermann}
\address{\newline Zimmermann \newline
Universit\'e de Picardie,
\newline D\'epartement de Math\'ematiques et LAMFA (UMR 7352 du CNRS),
\newline 33 rue St Leu,
\newline F-80039 Amiens Cedex 1,
\newline France}
\email{alexander.zimmermann@u-picardie.fr}

\dedicatory{Dedicated to Sri Wahyuni's $60^{th}$ birthday}

\newtheorem*{Theo2}{{Theorem}}
\newtheorem{Lemma1}{{Lemma}}
\newtheorem{Theo1}[Lemma1]{{Theorem}}
\newtheorem{Def1}[Lemma1]{{Definition}}
\newtheorem{Prop1}[Lemma1]{{Proposition}}
\newtheorem{Claim1}[Lemma1]{{Claim}}
\newtheorem{Rem1}[Lemma1]{{Remark}}
\newtheorem{Cor1}[Lemma1]{{Corollary}}
\newtheorem{Ex1}[Lemma1]{{Example}}
\newtheorem{Qu1}[Lemma1]{{Question}}

\newenvironment{Lemma}{\begin{Lemma1}}{\end{Lemma1}}
\newenvironment{Def}{\begin{Def1}\em}{\end{Def1}}
\newenvironment{Prop}{\begin{Prop1}}{\end{Prop1}}
\newenvironment{Rem}{\begin{Rem1}\rm}{\end{Rem1}}
\newenvironment{Theorem}{\begin{Theo1}}{\end{Theo1}}

\newenvironment{Cor}{\begin{Cor1}}{\end{Cor1}}

\newenvironment{Example}{\begin{Ex1}\em}{\end{Ex1}}

\newcommand{\uar}{\uparrow}
\newcommand{\dar}{\downarrow}
\newcommand{\lra}{\longrightarrow}

\newcommand{\ra}{\rightarrow}
\newcommand{\sdp}{\times\kern-.2em\vrule height1.1ex depth-.05ex}
\newcommand{\epi}{\lra \kern-.8em\ra}
\newcommand{\im}{\textup{im}}
\newcommand{\add}{\textup{add}}
\newcommand{\thick}{\textup{thick}}

\newcommand{\ul}{\underline}

\newcommand{\Z}{{\mathbb Z}}

 \normalbaselines

\date{\today}

\begin{document}

\begin{abstract}
For a finite group $G$ and an algebraically closed field $k$ of characteristic $p>0$
for any indecomposable finite dimensional $kG$-module $M$ with vertex $D$ and
a subgroup $H$ of $G$ containing $N_G(D)$ there is a unique indecomposable $kH$-module
$N$ of vertex $D$ being a direct summand of the restriction of $M$ to $H$.
This correspondence, called Green correspondence, was generalised by Auslander-Kleiner to
the situation of pairs of adjoint functors between additive categories. In the original
situation of group rings Carlson-Peng-Wheeler
proved that this correspondence is actually restriction of triangle functors between
triangulated quotient categories of the corresponding module categories. We review this theory
and show how we got a common generalisation of the approaches of Auslander-Kleiner
and Carlson-Peng-Wheeler, using Verdier localisations.
\end{abstract}

\maketitle

\section{Introduction}

Green correspondence is one of the most important tools in modular representation theory.
For a field $k$ of finite characteristic $p$ and a finite group $G$ we define
for any indecomposable $kG$-module $M$ its vertex and its source. Actually,
the vertex of $M$ is the smallest subgroup $D$ of $G$ such that $M$ is a direct factor of
$L\uar_D^G$ for some indecomposable $kD$-module $L$. This module $L$ is called the source
of $M$. Both $D$ and $L$ are essentially unique up to conjugacy, and it can be shown that
$D$ is always a $p$-subgroup of $G$. If $H$ contains the normaliser of $D$ in $G$, then
Green shows \cite{Green} that for any indecomposable $kG$-module $M$ with vertex $D$ there is a unique
indecomposable direct factor $f(M)$ of $M\dar^G_H$ with vertex $D$, its Green correspondent.
Carlson, Peng and Wheeler showed \cite{Carlson-Peng-Wheeler} much later that
this is actually the restriction of a functor of triangulated categories defined by some
well-known quotient construction imitating the construction of the standard stable category.
In a different direction Auslander and Kleiner showed \cite{AuslanderKleiner} that
Green correspondence actually works in a much more general setting of three additive
categories and pairs of adjoint functors between them. However, they did not mention
the triangulated category structure which is present underneath.
Further developments in this direction is given in \cite{Wang1}, which gives a
Green correspondence for derived categories.
In \cite{Wang2} a block version of the 
Green correspondence in the version from \cite{Carlson-Peng-Wheeler} is studied.
Most recently \cite{Wang3} studied Green correspondence for homotopy categories and derived categories of group rings
examining complexes of modules in a very detailed manner.
In \cite{Greentriangulated}
we gave a construction enlarging Auslander-Kleiner's approach to a triangulated category
situation. Specialising to the classical situation we get back the Carlson-Peng-Wheeler
theorem, generalising hence their result to a much vaster world. In the present note we
summarize and explain these approaches and discuss their links. We focus in particular
on the concept of $T$-relative projectivity, resp. $T$-relative injectivity
for a functor $T$ between triangulated categories. The classical
situation appears when we have a finite group $G$, a subgroup $H$ and
$T$ is the restriction functor from the (abelian)
category of finitely generated $kG$-modules to finitely generated $kH$-modules.

\paragraph{\bf Acknowledgement:}
I met Sri Wahynui for the first time in Autumn 2016, when I visited
Universitas Gadjah Mada in the context of
an exchange program initiated by Intan Muchtadi-Alamsyah and Indah Wijayanti. Sri
Wahyuni took responsibility for me during this
week and this was a most wonderful experience, mathematically, but
also and in particular personally. Intense
mathematical discussions over several hours with her and
various small groups of her students were followed by periods
where I got most impressively a very close
insight into the highly sophisticated
Javanese culture by her warm and generous personality.
I am most grateful to her for this experience which truly marks
one of the peek points in my mathematical career.
When Indah Wijayanti asked me some time later to speak on the
occasion of her 60th birthday I accepted right away without a second of hesitation.
The present paper deals with results around the talk I presented on the
occasion of her birthday party in Yogyakarta at 30 July 2019, and
I humbly wish to dedicate this paper to Sri Wahjuni.

\section{A Glimpse of Modular Representation Theory of Finite Groups; The Classical Case}

Let $G$ be a finite group and let $k$ be a field of characteristic $p>0$.
Denote by $kG$ the group ring of $G$ over the field $k$ (cf e.g. \cite[Chapter 1]{reptheobuch}).
Let $kG-mod$ be the category of finitely generated $kG$-modules.
Note that $kG-mod$, the category of finite-dimensional $kG$-modules
is a Krull-Schmidt category, that is, for every finite dimensional  $kG$-module
$M$ there are indecomposable $kG$-modules $M_1,\dots,M_s$ such that $M=M_1\oplus\dots\oplus M_s$,
and if there is another such set $N_1,\dots,N_t$ of indecomposable $kG$-modules
with $$N_1\oplus\dots \oplus N_t\simeq M\simeq M_1\oplus\dots\oplus M_s,$$
then $s=t$ and there is an element $\sigma$ in the symmetric group ${\mathfrak S}_t$ on $\{1,\dots,t\}$
such that $M_i\simeq N_{\sigma(i)}$ for all $i\in\{1,\dots,t\}$.

Note moreover that for any algebra $A$ a $d$-dimensional $A$-module is the same as a $k$-algebra homomorphism
$$A\stackrel{\mu}{\lra}End_k(k^d). $$
Hence, for any $k$-algebra homomorphism $B\stackrel{\alpha}{\lra}A$ and
any $d$-dimensional $A$-module $M$ we obtain a $d$-dimensional $B$-module
$res^A_\alpha(M)$, or $res_B^A$ for short if $\alpha$ is evident from the context,
given by the algebra homomorphism $\mu\circ\alpha$.

For a subgroup $H$ of $G$ there is an algebra homomorphism $\iota:kH\lra kG$
given by the inclusion map of $H$ in $G$.
By the above any $kG$-module $M$ induces a $kH$-module $res^{kG}_{\iota}(M)$, denoted also by
$M\dar^G_H$. We call this the restriction of $M$ to $H$.
It is clear that $res^G_H$ is a functor $kG-mod\lra kH-mod$ since
if a map $M_1\lra M_2$ is $kG$-linear, then it is
trivially $kH$-linear. Moreover,
$res^G_H$ admits a left adjoint $ind_H^G:kH-mod\lra kG-mod$ in the sense that
for any $kG$-module $M$ and any $kH$-module $N$ we have
$$Hom_{kH}(N,res^G_H(M))\simeq Hom_{kG}(ind_H^G(N),M)$$
and this isomorphism is functorial in $M$ and $N$. Actually $kG$ is a $kG-kH$-bimodule, and then
$$res^G_H(M)=Hom_{kG}(kG,M).$$
Therefore $ind_H^G(N)=kG\otimes_{kH}N$
and the adjointness is nothing else than an incident of the usual
Hom-tensor adjunction.
A more detailed analysis (cf \cite[Chapter 1]{reptheobuch})
shows that $ind_H^G$ is also right adjoint to
$res^G_H$. To be slightly more precise,
the right adjoint of $res^G_H$ is coinduction, and if $G$ is a finite group, then coinduction
equals induction.

If $M$ is an indecomposable $kG$-module, then $M\dar^G_H$ does not need to be
an indecomposable $kH$-module. An easy example is $H=\{1\}$, then $M\dar^G_H$ is indecomposable
if and only if $M$ is one-dimensional.

\subsection{Relative Projectivity and Vertices; Classical Case}

A key notion in the above context is relative projectivity. This is classical, and
the definition we use, follows Hoch\-schild's work \cite{Hochschild}. We shall study his concept in more detail in Section~\ref{Hochschilddefinition} below.

\begin{Def}\label{relativeprojectiveclassicalDef}
Let $G$ be a finite group, let $H$ be a subgroup of $G$, and let
$M$ be a finite dimensional $kG$-module. Then $M$ is relatively $H$-projective if for all
$kG$-modules $N$ and all epimorphisms $\alpha:N\lra M$ such that the induced morphism
$res(\alpha):res^G_H(N)\lra res^G_H(M)$ is split, then also $\alpha$ is split.
\end{Def}

The most interesting result is now the following statement, known as Higman's lemma.

\begin{Lemma}\label{Higmanslemma}
Let $G$ be a finite group, let $k$ be a field, and let $H$ be a subgroup of $G$.
Let $M$ be a finite-dimensional indecomposable $kG$-module.
Then the following statements are equivalent.
\begin{itemize}
\item $M$ is
relatively $H$-projective
\item  $M$ is a direct summand of $ind_H^G(res^G_H(M))$
\item There is a finite-dimensional indecomposable $kH$-module $L$
such that $M$ is a direct summand of $ind_H^G(L)$.
\end{itemize}
\end{Lemma}

A first observation is that $M$ is a direct summand of $ind_1^G(res^G_1(M))$ if and
only if $M$ is projective. Actually, any module of the form $ind_1^G(L)$ is free of rank $\dim_kL$,
and the inverse implication follows by Frobenius reciprocity and the fact that
vector spaces are always free. Now, the question arises, for a given indecomposable
and not necessarily projective
$kG$-module $M$, what are minimal subgroups $H$
such that $M$ is a direct factor of $ind_H^G(res^G_H(M))$.

Using Lemma~\ref{Higmanslemma} and Mackey's decomposition, i.e. the decomposition
of $kG$ into indecomposable $kH-kH$-bimodules,
the following consequence is immediate.

\begin{Theorem}\label{vertexexistence}
Let $G$ be a finite group, let $k$ be a field of
characteristic $p>0$, and let $H$ be a subgroup of $G$.
Then for any indecomposable $kG$-module $M$ the set
of minimal elements in the partial ordered
set of subgroups $H$ such that $M$ is relatively
$H$-projective forms a $G$-conjugacy class of
$p$-subgroups of $G$.
\end{Theorem}

Let now $k$ be a field of characteristic $p>0$.
For an indecomposable $kG$-module $M$ we know by Theorem~\ref{vertexexistence}
that if $M$ is relatively $D$-projective and if $D$ is minimal with this property, then
$D$ is a $p$-subgroup of $G$, and two of these subgroups are conjugate in $G$.
We call such a subgroup $D$ a {\em vertex of $M$}. Moreover Higman's Lemma~\ref{Higmanslemma}
shows that for an indecomposable $kG$-module $M$ with vertex $D$
there is an indecomposable $kD$-module $L$ such that $M$ is a direct summand of $M$.
We call $L$ a {\em source of $M$.} Arguments using essentially the Krull-Schmidt theorem shows that also
a source is basically unique.

\subsection{Green Correspondence; The Classical Case}

One of the most important and basic tools in representation theory of
finite groups, namely Green correspondence, is attached to these concepts.
As usual for a group $G$, a subgroup $H$ of $G$ and an element $g\in G$ we denote
$^gH:=\{ghg^{-1}\;|\;h\in H\}$. The following Theorem~\ref{Greencorrespondenceclassical} is called Green correspondence.

\begin{Theorem}(Green~\cite{Green})\label{Greencorrespondenceclassical}
Let $G$ be a finite group, let $k$ be a field of characteristic
$p>0$, and let $D$ be a $p$-subgroup of $G$.
Let $H$ be a subgroup of $G$ containing $N_G(D)$. Put
$${\mathfrak X}:=\{S\leq D\cap {}^gD\;|\;g\in G\setminus H\}$$
$${\mathfrak Y}:=\{S\leq H\cap {}^gD\;|\;g\in G\setminus H\}$$
Then
\begin{itemize}
\item for any indecomposable $kG$-module $M$ with vertex $D$
there is a unique indecomposable direct summand $f(M)$ of $res^G_H(M)$ with vertex $D$.
All other indecomposable direct summands of $res^G_H(M)$ have vertex in $\mathfrak Y$.
\item for any indecomposable $kH$-module $N$ with vertex $D$
there is a unique indecomposable direct summand $g(N)$ of $ind^G_H(N)$ with vertex $D$.
All other indecomposable direct summands of $ind^G_H(N)$ have vertex in $\mathfrak X$.
\item $fg=id$ and $gf=id$.
\end{itemize}
\end{Theorem}

\begin{Example}\label{ExampleTIdef}
Consider the special case of $D$ being cyclic of order $p$. Then
$D$ has only one proper subgroup, namely the trivial group. Then, for any
indecomposable $kG$-module with vertex $D$,
by Theorem~\ref{Greencorrespondenceclassical} $res^G_H(M)$ has a unique
non projective direct summand $f(M)$, and all other indecomposable direct factors are projective.
A similar special case is given for $G$ being a group with trivial intersection property of
Sylow subgroups $D$, that is $\{{}^gD\cap D\:|\;g\in G\}\subseteq\{1,D\}$. There,
the same statement holds.
\end{Example}

This observation motivates the following construction, namely the
well-established tool in representation theory, the  stable category.
For a finite dimensional algebra $A$ let $A-\ul{mod}$ be the category whose objects are $A$-modules.
For any two $A$-modules let $PHom_A(M,N)$ be the vector space of $A$-module
homomorphisms $M\lra N$ which factor through a projective $A$-module.
Then, the morphisms from $M$ to $N$ in the stable category are elements in
$\ul{Hom}_A(M,N):=Hom_A(M,N)/PHom_A(M,N)$. Composition of morphisms in the stable category
is given by composition of morphisms of $A$-modules, and it is easy to see that this gives
a well-defined category. Projective $A$-modules are isomorphic to $0$, since the endomorphism
algebra in the stable category of a projective module is $0$. Moreover, most interestingly,
if $A$ is self-injective, then $A-\ul{mod}$ is a triangulated category (cf e.g.~\cite[Chapter 3]{reptheobuch}).

\begin{Example}\label{ExampleTIstablecat}
We come back to Example~\ref{ExampleTIdef}. Restriction and induction of projective modules are projective. Hence
restriction and induction induce functors on the level of the stable categories.
$$ind_H^G:kH-\ul{mod}\lra kG-\ul{mod}$$
$$res_H^G:kG-\ul{mod}\lra kH-\ul{mod}$$
Consider Green correspondence for $D$ a cyclic group of order $p$, or
if $G$ is a finite group with trivial intersection Sylow $p$-subgroups and $D$ is a Sylow subgroup. Put $H$ a
subgroup of $G$ containing $N_G(D)$.
Then Theorem~\ref{Greencorrespondenceclassical} shows that in these cases $res^G_H$ and $ind_H^G$ induce
equivalences of categories of $kG$-modules with vertex $D$ and $kH$-modules with vertex $D$.
Observe that neither $ind_H^G$ nor $res^G_H$ are equivalences of categories, but they are
equivalences of the additive subcategories generated by indecomposable modules of vertex $D$.
Note that these subcategories are not triangulated subcategories of the stable category in general.
\end{Example}

\subsection{Green Correspondence is the Trace of Triangle Functors}

\label{section2}

In \cite{Carlson-Peng-Wheeler} Carlson-Peng-Wheeler showed that Green correspondence is actually
the restriction of a triangle functor between certain triangulated subcategories of the corresponding module categories over the relevant groups.

More precisely, for an additive category $\mathcal A$ and an additive subcategory $\mathcal S$
denote by ${\mathcal A}/{\mathcal S}$ the category with the same objects as $\mathcal A$. A morphism $f$ in ${\mathcal A}(X,Y)$ is said to be in ${\mathcal A}^{\mathcal S}(X,Y)$ if there
an object $Z$ of $\mathcal S$ and morphisms $g\in{\mathcal A}(X,Z)$ and $h\in{\mathcal A}(Z,Y)$
such that $f=h\circ g$. Then put
$$({\mathcal A}/{\mathcal S})(X,Y):={\mathcal A}(X,Y)/{\mathcal A}^{\mathcal S}(X,Y)$$
and composition is given by composition of representatives of classes. It is clear that this is well-defined, such as in the remarks following Example~\ref{ExampleTIdef}.

Carlson-Peng-Wheeler define for a fixed finite dimensional $kG$-module $W$ a
finite dimensional module $M$ to be $W$-projective if
$M$ is a direct factor of $W\otimes_kW^*\otimes M$, where as usual $-^*$ denotes the
$k$-linear dual. Then, they show

\begin{Prop}\cite[Section 6]{Carlson-Peng-Wheeler} Let $k$ be a field and let $G$ be a finite group.
Let $W$ be a finitely generated $kG$-module and let $kG-mod^W$ be the full subcategory of $W$-projective $kG$-modules. Then $kG-mod/kG-mod^W$ carries the structure of
a triangulated category.
\end{Prop}

As a consequence they show that Green correspondence is the restriction of triangle
functors of quotient categories of $kG-mod$
respectively $kH-mod$ for appropriate choices of $W$. More precisely

\begin{Theorem} \cite{Carlson-Peng-Wheeler} \label{GreenCarlson-Peng-Wheeler}
Let $k$ be a field of characteristic $p>0$, let $G$ be a finite group, let $D$ be a $p$-subgroup of $G$, and let
$H$ be a subgroup of $G$ containing $N_G(D)$. For any group $\Gamma$ denote by $k$ the trivial $k\Gamma$-module.
Put
$${\mathfrak X}:=\{S\leq D\cap {}^gD\;|\;g\in G\setminus H\}$$
$${\mathfrak Y}:=\{S\leq H\cap {}^gD\;|\;g\in G\setminus H\}$$
and $W_G:=\bigoplus_{X\in{\mathcal X}}k\uar_X^G$ and $W_H:=\bigoplus_{Y\in{\mathcal Y}}k\uar_Y^H$.
Then restriction $\dar^G_H$ induces a triangle functor $$kG-mod/{kG-mod}^{W_G}\lra kH-mod/{kH-mod}^{W_H}$$
and induction $uar_H^G$ induces a triangle functor $$kH-mod/{kH-mod}^{W_H}\lra kG-mod/{kG-mod}^{W_G}.$$
These functors restrict to equivalences between the full additive subcategories generated by indecomposable
modules of vertex $D$.
\end{Theorem}

Note that this version of Green correspondence is not an equivalence of triangulated categories, but that
the Green correspondence is the restriction of triangle functors between triangulated categories.

\section{Notions of Relative Projectivity}

We observe that in order to generalize Green correspondence to a categorical setting we first need to
generalise and formulate the notion of relative projectivity.
The classical case provides two such settings.

\subsection{Recall Hochschild's concept}

\label{Hochschilddefinition}

We recall Hochschild's approach \cite[Section 1]{Hochschild} to
relative projectivity (respectively injectivity).
He considers an algebra $R$ and a subalgebra $S$ and says that an exact sequence
$$\xymatrix{\cdots\ar[r]&M_{i+2}\ar[r]^{t_{i+2}}&M_{i+1}\ar[r]^{t_{i+1}}&M_i\ar[r]^{t_i}& M_{i-1}\ar[r]&\dots}$$
of $R$-module homomorphisms is said to be
$(R,S)$-exact if the kernel of $t_i$ is a direct summand as $S$-modules of $M_i$ for all $i$.
Equivalently the sequence is $(R,S)$-exact if
$t_i\circ t_{i+1}=0$ for all $i\in\Z$ and in addition
there are $S$-module homomorphisms $h_i:M_i\ra M_{i+1}$ such that
$t_{i+1}\circ h_i+h_{i-1}\circ t_i=id_{M_i}$ for all $i\in\Z$.
Hochschild continues that an $R$-module $A$ is called $(R,S)$-injective
if for every $(R,S)$-exact sequence
$$\xymatrix{0\ar[r]&U\ar[r]^p&V\ar[r]^q&W\ar[r]&0}$$
and every $R$-module homomorphism $h:U\ra A$ there is an $R$-module homomorphism
$h':V\ra A$ with $h=h'\circ p$. Dually, an $R$-module $A$ is $(R,S)$-projective, if for each $R$-module homomorphism $g:A\ra W$ there is an $R$-module homomorphism
$g':A\ra V$ such that $g=q\circ g'$.

What precisely are short $(R,S)$-exact sequences?

\begin{Lemma}
An exact sequence
$$\xymatrix{0\ar[r]&U\ar[r]^p&V\ar[r]^q&W\ar[r]&0}$$
of $R$-modules is $(R,S)$-exact if, and only if, it splits when considered as
sequence of $S$-modules.
\end{Lemma}

\begin{proof}
If the sequence is $R-S$-exact, then there is an $S$-homotopy $h_1:W\ra V$ such that $h_1\circ q=id_W$, whence the sequence splits. Conversely if the restriction of the sequence splits, then
by definition the kernel of $q$ is an $S$-direct summand of $V$, namely $p(U)$.
\end{proof}

\begin{Lemma}
 $A$ is $(R,S)$-injective
if, and only if, each $(R,S)$-exact sequence $$\xymatrix{0\ar[r]&A\ar[r]^p&V\ar[r]^q&W\ar[r]&0}$$
splits as sequence of $R$-modules. Similarly,
$A$ is $(R,S)$-projective
if, and only if, each $(R,S)$-exact sequence $$\xymatrix{0\ar[r]&U\ar[r]^p&V\ar[r]^q&A\ar[r]&0}$$
splits as sequence of $R$-modules.
\end{Lemma}

\begin{proof} This is analogous to the usual argument in homological algebra.
Indeed, if $A$ is $(R,S)$-projective, and
$$\xymatrix{0\ar[r]&U\ar[r]^p&V\ar[r]^q&A\ar[r]&0}$$
is an $(R,S)$-exact sequence, then by definition, for $id_A$ there is an $R$-module homomorphism $s:A\ra V$ with $id_A=q\circ s$. This is tantamount to say that the sequence splits.
Conversely, suppose that each $(R,S)$-exact sequence
$$\xymatrix{0\ar[r]&U'\ar[r]^{p'}&V'\ar[r]^{q'}&A\ar[r]&0}$$
splits. Let
$$\xymatrix{0\ar[r]&U\ar[r]^p&V\ar[r]^q&W\ar[r]&0}$$
be an $(R,S)$-exact sequence and consider an $R$-module homomorphism $g:A\ra W$.
Then form the pullback of the sequence
$$\xymatrix{0\ar[r]&U\ar[r]^p&V\ar[r]^q&W\ar[r]&0}$$ along $g$ to get a commutative diagram
$$
\xymatrix{0\ar[r]&U\ar[r]^p&V\ar[r]^q&W\ar[r]&0\\
0\ar[r]&U\ar[r]^{p'}\ar@{=}[u]&V'\ar[r]^{q'}\ar[u]^x&A\ar[r]\ar[u]^g&0}
$$
and the bottom sequence is again $(R,S)$-exact, and hence splits, by hypothesis.
Therefore there is an $R$-module homomorphism $s:A\ra V'$ with $q'\circ s=id_A$.
Hence $q\circ x\circ s=g\circ q'\circ s=g$ as required.

The case of $(R,S)$-injective is dual.
\end{proof}

\begin{Cor}
An $R$-module $A$ is $(R,S)$-projective if and only if each short exact sequence
of $R$-modules
$$\xymatrix{0\ar[r]&U\ar[r]^p&V\ar[r]^q&A\ar[r]&0}$$
which is known to split as sequence of $S$-modules, is also split as sequence of $R$-modules.

An $R$-module $A$ is $(R,S)$-projective if and only if each short exact sequence
of $R$-modules
$$\xymatrix{0\ar[r]&A\ar[r]^p&V\ar[r]^q&W\ar[r]&0}$$
which is known to split as sequence of $S$-modules, is also split as sequence of $R$-modules.

In more modern terms, denoting by $res^R_S:R-Mod\lra S-Mod$ the restriction functor,
then $A$ is $(R,S)$-projective if and only if the induced functor
$$Ext^1_R(A,-)\lra Ext^1_S(res^R_S(A),res^R_S(-))$$
is a monomorphism in the functor category $R-Mod\lra \Z-Mod$.
An $R$-module $A$ is $(R,S)$-injective if and only if the induced functor
$$Ext^1_R(-,A)\lra Ext^1_S(res^R_S(-),res^R_S(A))$$
is a monomorphism in the functor category $R-Mod\lra \Z-Mod$.
\end{Cor}

\subsection{The new concept for triangulated categories}

In this subsection we shall give an alternative proof for the results
in \cite[Section 2]{Greentriangulated}.
Let $\mathcal S$ and $\mathcal T$ be exact categories, and let
$S:{\mathcal T}\ra{\mathcal S}$ be an exact functor.
Denote by $Ext^1_{\mathcal T}(X,Y)$ be the set of equivalence classes of short exact sequences
$$0\lra Y\stackrel{\iota}{\lra} E\stackrel{\pi}\lra X\lra 0,$$
where as usual two such sequences
$$0\lra Y\stackrel{\iota_1}{\lra} E_1\stackrel{\pi_1}\lra X\lra 0,$$
and
$$0\lra Y\stackrel{\iota_2}{\lra} E_2\stackrel{\pi_2}\lra X\lra 0,$$
are equivalent if and only if there is a homomorphism $E_1\lra E_2$ making the diagram
$$
\xymatrix{
0\ar[r]& Y\ar[r]^{\iota_1}\ar@{=}[d]& E_1\ar[r]^{\pi_1}\ar[d]& X\ar[r]\ar@{=}[d]& 0\\
0\ar[r]& Y\ar[r]^{\iota_2}& E_2\ar[r]^{\pi_2}& X\ar[r]& 0
}
$$
commutative. Note that by \cite[Corollary 3.2]{Buehler} any such homomorphism
$E_1\lra E_2$ is an isomorphism.
This way $S$ induces a morphism
$$S:Ext_{\mathcal T}^1(X,Y)\lra Ext_{\mathcal S}^1(SX,SY)$$
Now, if ${\mathcal S}=A-mod$ for some algebra $A$, then denoting by $D^b(A)$
the derived category of bounded complex of finitely generated $A$-modules
(cf e.g. \cite[Chapter 3]{reptheobuch}),
$$Ext^1_{\mathcal S}(X,Y)=Ext^1_{A}(X,Y)=Hom_{D^b(A)}(X,Y[1]).$$
Consider Definition~\ref{relativeprojectiveclassicalDef}. A module $M$ is
relative $H$-projective if for all modules $X$
the map $$res^G_H:Ext^1_{kG}(M,X)\lra Ext^1_{kH}(res^G_H(M),res^G_H(X))$$
is injective. In other words,
a module $M$ is
relative $H$-projective if the natural transformation
$$res^G_H:Ext^1_{kG}(M,-)\lra Ext^1_{kH}(res^G_H(M),res^G_H(-))$$
is a monomorphism in the category of functors $kG-mod\lra k-mod$.

If we want to enlarge the notion of relative projectivity to the derived category
we first observe that $res^G_H$ is a triangle functor $S:{\mathcal T}\lra{\mathcal S}$
from the triangulated category ${\mathcal T}:=D^b(kG)$ to the triangulated category
${\mathcal S}=D^b(kH)$, and the notion should be formulated
for triangle functors. Then we see that an object $M$ in $D^b(A)$ is
$S$-relative projective if and only if the natural transformation
$$S:Hom_{\mathcal T}(M,-[1])\lra Hom_{\mathcal S}(SM,S-[1])$$
is a monomorphism of functors ${\mathcal T}\lra k-mod$. Since this functor
can be evaluated on all objects of ${\mathcal T}$, and since $[1]$ is an
auto-equivalence of $\mathcal T$, we can just omit $[1]$ in the above formula.

\begin{Def} \cite{Greentriangulated}
Let $\mathcal T$ and $\mathcal S$ be triangulated categories, and let $S:{\mathcal T}\lra{\mathcal S}$
be a triangle functor. Then
\begin{itemize}
\item
an object $M$ of $\mathcal T$ is $S$-relatively projective if the natural transformation
$$S:Hom_{\mathcal T}(M,-)\lra Hom_{\mathcal S}(SM,S-)$$
is a monomorphism in the functor category ${\mathcal T}\lra k-mod$.
\item
an object $N$ of $\mathcal T$ is $S$-relatively injective if the natural transformation
$$S:Hom_{\mathcal T}(-,N)\lra Hom_{\mathcal S}(S-,SN)$$
is a monomorphism in the functor category ${\mathcal T}^{op}\lra k-mod$.
\end{itemize}
\end{Def}

Of course, this is quite a large generalisation of the classical notion of relative projectivity.
The classical case is found for self-injective algebras, such as group algebras, for the stable
category rather than the module category. Indeed, we omitted by purpose the shift of degree by $1$.
But then $Ext^0$ is part of our study, and this should be the stable homomorphisms, and not
just the ordinary homomorphisms.

How what the alternative definition of relative projectivity using Green's definition coming
from Higman's lemma. Is there some link, or a triangulated version of Higman's lemma?
Most astonishing, this is true, at least in the correct setting. For a subcategory
$\mathcal S$ of an additive category $\mathcal C$ denote by $\add({\mathcal S})$
the additive closure of $\mathcal S$ in $\mathcal C$.

\begin{Rem}
We emphasize that we denote by $\epsilon$ the unit of an adjunction and by $\eta$ the counit of an adjunction. Many papers use the inverse convention, and we alert the reader to pay attention to this fact.
\end{Rem}

\begin{Lemma}\label{relatively-t-projective-for-adjoints} \cite{Greentriangulated}
Let ${\mathcal T}$ and ${\mathcal S}$ be $k$-linear categories and let $T:{\mathcal S}\lra {\mathcal T}$
be a $k$-linear functor. Suppose that $T$ has a left (respectively right) adjoint $S$, and denote by
$$\varphi_{X,Y}:{\mathcal T}(X,TY)\stackrel{\simeq}\lra {\mathcal S}(SX,Y)$$
the adjunction isomorphism.
Then any object in $\add(\im(S))$ is $T$-relative projective (respective injective).
\end{Lemma}

\begin{proof} Let $\epsilon:id_{\mathcal T}\ra TS$ be the unit of the adjunction.
By \cite[IV Theorem 1.(i)]{Maclane}  for any $f\in{\mathcal S}(SX,Y)$ we have
$$\varphi_{X,Y}^{-1}(f)=T(f)\circ\epsilon_X.$$
We first suppose $Q=SQ'$ for some object $Q'$ of $\mathcal T$.
Consider the following diagram
$$
\xymatrix{{\mathcal S}(Q,-)\ar[r]^-{T_Q}\ar@{=}[d]&{\mathcal T}(TQ,T-)\ar@{=}[d]\\
{\mathcal S}(SQ',-)\ar[r]^-{T_{SQ'}}&{\mathcal T}(TSQ',T-)\ar[d]^{{\mathcal T}(\epsilon_{Q'},-)}\\
{\mathcal T}(Q',T-)\ar[u]^{\varphi_{Q',-}}\ar@{-->}[r]^-\lambda&{\mathcal T}(Q',T-)}
$$
We define $\lambda:={\mathcal T}(\epsilon_{Q'},-)\circ T_{SQ'}\circ\varphi_{Q'}$,
as indicated in the above diagram.
Hence
$$\lambda(f)=T_{SQ'}(\varphi_{Q',-}(f))\circ\epsilon_{Q'}=\varphi_{Q',-}^{-1}(\varphi_{Q',-}(f))=f.$$
Since $\varphi_{Q',-}$ is an isomorphism, $T_{Q}$ is split mono.
Let now $Q\in \add(\im(S))$ and $Q$ is a direct factor of $SQ'$. We simply use
that the above argument is still valid for direct factors of $SQ'$.

The case of relative injective is done analogously, using the counit instead of the unit.
\end{proof}

\begin{Prop}\label{relativelyprojectivecharacterisation} \cite{Greentriangulated}
Let ${\mathcal T}$ and ${\mathcal S}$ be triangulated categories and let $T:{\mathcal S}\lra {\mathcal T}$
be a triangle functor. Suppose that $T$ has a left (respectively right) adjoint $S$.
Then an object $Q$ is $T$-relative projective (respectively injective) if and only if $Q$ is in $\add(\im(S))$.
\end{Prop}

\begin{proof}
By Lemma~\ref{relatively-t-projective-for-adjoints} we see that any object in $\add(\im(S))$
is $T$-relative projective.

\medskip

Suppose now that $Q$ is $T$-relative projective.
Let $\eta:ST\lra id_{\mathcal S}$ be the counit of the adjunction. We shall show that
$\eta_Q$ is a split epimorphism. This then shows that $Q$ is in $\add(\im(S))$.

Again by \cite[IV Theorem 1]{Maclane} the composition
$$\xymatrix{T\ar[r]^{\epsilon{T}}&TST\ar[r]^{T\eta}&T}$$
is the identity, and so $T\eta$ is a split epimorphism. Let
$$\xymatrix{\widetilde Q\ar[r]&STQ\ar[r]^{\eta_Q}&Q\ar[r]^-\nu&\widetilde Q[1]}$$
be a distinguished triangle. Then
$$\xymatrix{T\widetilde Q\ar[r]&TSTQ\ar[r]^{T\eta_Q}&TQ\ar[r]^-{T\nu}&T\widetilde Q[1]}$$
is a distinguished triangle as well. Since $T\eta_Q$ is a split epimorphism,
$T\nu=0$. Now, $Q$ is $T$-relative projective, and hence by hypothesis
$$\xymatrix{{\mathcal S}(Q,-)\ar[r]^-{T_Q}&{\mathcal T}(TQ,T-)}$$
is injective. Evaluate this on $\widetilde Q[1]$ and obtain that $\nu=0$. Hence, $\eta_Q$ is
split epimorphism, which shows the statement.

The proof of the case of $T$-relative injective objects is done completely analogously
using the sequence
$$\xymatrix{T\ar[r]^{T\widetilde\epsilon}&TST\ar[r]^{\widetilde{\eta}{T}}&T}$$
for the unit $\widetilde \epsilon$ and the counit $\widetilde\eta$ of the adjunction $(T,S)$.
\end{proof}

\begin{Rem}
Recall that $Q\mapsto \widetilde Q$ is {\em not} a functor.
\end{Rem}

\begin{Cor}\label{counitsplitsmeansrelativelyprojective}
Let ${\mathcal T}$ and ${\mathcal S}$ be triangulated categories and let $T:{\mathcal S}\lra {\mathcal T}$
be a triangle functor. Suppose that $T$ has a left (respectively right) adjoint $S$, and let $\eta:ST\lra\textup{id}$ be the
counit (respectively $\widetilde\epsilon:\textup{id}\lra ST$ the unit) of the adjunction.
Then $Q$ is $T$-relative projective (respectively injective) if and only if $\eta_Q$ is a split epimorphism
(respectively $\widetilde\epsilon_Q$ is a split monomorphism).
\end{Cor}

\begin{proof}
Suppose that $Q$ is $T$-relative projective. By Proposition~\ref{relativelyprojectivecharacterisation}
we see that $Q$ is in $\add(\im S)$. By \cite[IV Theorem 1]{Maclane} $\eta_{SR}$ is a split epimorphism
for any object $R$ of $\mathcal T$. Hence $\eta_Q$ is a split epimorphism.
%Let
%$$\xymatrix{WQ\ar[r]&STQ\ar[r]^{\eta_Q}&Q\ar[r]^-\nu&WQ[1]}$$
%be a distinguished triangle. The second part of the proof of
%Proposition~\ref{relativelyprojectivecharacterisation} shows that $\nu=0$ and hence $\eta_Q$ is
%a split epimorphism.

If $\eta_Q$ is a split epimorphism, then $Q$ is in $\add(\im(S))$, and by
Proposition~\ref{relativelyprojectivecharacterisation} this implies that $Q$ is $T$-relatively
projective.
\end{proof}

We summarize the situation to an analogue of Higman's lemma for pairs of adjoint functors
between triangulated categories.

\begin{Theorem} \cite[Proposition 2.7]{Greentriangulated}
\label{HigmansLemmaForAdjointFunctors}
Let ${\mathcal T}$ and ${\mathcal S}$ be triangulated categories and let $T:{\mathcal S}\lra {\mathcal T}$
be a triangle functor. Suppose that $T$ has a left (respectively right) adjoint $S$. Let $M$ be an
indecomposable object of $\mathcal T$. Then the following
are equivalent:
\begin{enumerate}
\item $M$ is $T$-relatively projective (respectively injective).
\item $M$ is in $\add(\im S)$.
\item $M$ is a direct factor of some $S(L)$ for some $L$ in $\mathcal S$.
\item $M$ is a direct factor of $ST(M)$.
\end{enumerate}
\end{Theorem}

\begin{proof}
$(1)\Leftrightarrow (2)$ is Proposition~\ref{relativelyprojectivecharacterisation}.

$(2)\Leftrightarrow (3)$ is the definition of $\add(\im S)$.

$(3)\Rightarrow (4)$ is trivial.

$(4)\Rightarrow (1)$ is Corollary~\ref{counitsplitsmeansrelativelyprojective}.
\end{proof}

\begin{Rem}
Note that Corollary~\ref{counitsplitsmeansrelativelyprojective} generalises
\cite[Proposition 2.1.6, Proposition 2.1.8]{reptheobuch} to this more general situation.
\end{Rem}

Note that condition $(2)$ indicates that in case $T$ has a left adjoint $S_\ell$,
a right adjoint $S_r$, and in case $S_\ell=S_r$, then $T$-relatively
injective and $T$-relatively projective is just the same property. The situation
also occurs under the weaker condition $\add(\im S_\ell)=\add(\im S_r)$.

\begin{Rem}
In \cite{Greentriangulated} we restricted the notion of $T$-relative projectivity
respectively $T$-relative injectivity
to the case of functors $T$ having a left respectively right adjoint. This
is caused by the
fact that in \cite{Greentriangulated} we are guided there by the approach of Beligiannis-Marmaridis
\cite{BeligiannisMarmaridis}, whereas here we rather use Hochschild's approach.
By Theorem~\ref{HigmansLemmaForAdjointFunctors} and the corresponding statement
in \cite{Greentriangulated} both definitions give the same result
if $T$ has a left and a right adjoint.
If $T$ has a left and right adjoint then Brou\'e showed in \cite{BroueHigman}
a  a slightly different version of 
Theorem~\ref{HigmansLemmaForAdjointFunctors} by completely different methods. 
\end{Rem}

\section{Auslander Kleiner's Version of Green Correspondence}

Auslander and Kleiner proposed in \cite{AuslanderKleiner}
a version of Green correspondence which worked for pairs of adjoint functors between
abelian categories. They observed that the arguments of classical
Green correspondence are essentially disguised arguments on pairs of adjoint functors.

We start with some notations.
\begin{itemize}
\item
Recall from Section~\ref{section2} the following notation.
Let $\mathcal U$ be an additive category, and let $\mathcal V$ be a full subcategory of $\mathcal U$.
Then we denote by ${\mathcal U}/{\mathcal V}$ the category whose objects are the same as
the objects of $\mathcal U$, and for any two objects $X$ and $Y$ of
${\mathcal U}/{\mathcal V}$ the morphisms from $X$ to $Y$ in ${\mathcal U}/{\mathcal V}$
are ${\mathcal U}(X,Y)/I_{\mathcal U}^{\mathcal V}(X,Y)$, where
$I_{\mathcal U}^{\mathcal V}(X,Y)$ is the subset of ${\mathcal U}(X,Y)$ given by
those $f\in{\mathcal U}(X,Y)$ such that there is an object $Z\in{\mathcal V}$ and
$h\in {\mathcal U}(X,Z)$, $g\in{\mathcal U}(Z,Y)$ such that $f=g\circ h$.
\item
Let ${\mathcal A}$ and $\mathcal B$ be additive categories, and let $F:{\mathcal A}\lra{\mathcal B}$
be a functor. Then for any full subcategory $\mathcal C$ of $\mathcal B$ let
$F^{-1}({\mathcal C}$ be the full subcategory of $\mathcal A$ generated by the objects $X$ of $\mathcal A$
such that $F(X)$ is a direct summand of an object of $\mathcal C$.
\item
If ${\mathcal S}$ and $\mathcal R$ are subcategories of a Krull-Schmidt category $\mathcal W$, then
${\mathcal R}-{\mathcal S}$ denotes the full subcategory of ${\mathcal R}$ consisting of those
objects $X$ of $\mathcal R$ such that no direct factor of $X$ is an object of $\mathcal S$.
\end{itemize}

Let ${\mathcal D}$, ${\mathcal H}$, ${\mathcal G}$ be three additive categories
$$
\xymatrix{
{\mathcal D}\ar@/^/[r]^{S'}&{\mathcal H}\ar@/^/[l]^{T'}\ar@/^/[r]^{S}&{\mathcal G}\ar@/^/[l]^T
}
$$
such that $(S,T)$ and $(S',T')$ are adjoint pairs. Let $\epsilon:id_{\mathcal H}\lra TS$
be the unit of the adjunction $(S,T)$.

We assume throughout the rest of the section that there is an
endofunctor $U$ of $\mathcal H$ such that $TS=1_{\mathcal H}\oplus U$, denote by
$p_1:TS\lra 1_{\mathcal H}$
the projection, and suppose that $p_1\circ\epsilon$ is an isomorphism.
Note that we use both of the notations $id_{\mathcal C}$ and $1_{\mathcal C}$ for the identity functor
on the category $\mathcal C$.

\begin{Theorem} \cite[Theorem 1.10]{AuslanderKleiner} \label{Greenforadjoints}
Assume that there is an
endofunctor $U$ of $\mathcal H$ such that $TS=1_{\mathcal H}\oplus U$, denote by
$p_1:TS\lra 1_{\mathcal H}$
the projection, and suppose that $p_1\circ\epsilon$ is an isomorphism.

Let $\mathcal Y$ be a subcategory of $\mathcal H$ and let
${\mathcal Z}:=(US')^{-1}({\mathcal Y})$.
Then the following two conditions $(\dagger)$ are equivalent.
\begin{itemize}
\item
Each object of $S'T'{\mathcal Y}$ is a direct factor of an object of
$\mathcal Y$ and of an object of $U^{-1}({\mathcal Y})$.
\item Each object of $TSS'T'{\mathcal Y}$ is a direct factor of an object of
$\mathcal Y$.
\end{itemize}
Suppose that the above conditions $(\dagger)$ hold for $\mathcal Y$. Then
\begin{enumerate}
\item $S$ and $T$ induce functors
$${\mathcal H}/S'T'{\mathcal Y}\stackrel{S}{\lra} {\mathcal G}/SS'T'{\mathcal Y}\mbox{ and }
{\mathcal G}/SS'T'{\mathcal Y}\stackrel{T}{\lra}{\mathcal H}/{\mathcal Y}$$
\item For any object $L$ of $\mathcal D$ and any object $B$ of $U^{-1}({\mathcal Y})$ the functor
$S$ induces an isomorphism
$${\mathcal H}/S'T'{\mathcal Y}(S'L,B)\lra {\mathcal G}/SS'T'{\mathcal Y}(SS'L,SB)$$
\item
For any object $L$ of $(US')^{-1}{\mathcal Y}$ and any object $A$ of ${\mathcal G}$ the functor
$T$ induces an isomorphism
$${\mathcal G}/SS'T'{\mathcal Y}(SS'L,B)\lra {\mathcal H}/{\mathcal Y}(TSS'L,TA)$$
\item\label{item4}
The restrictions of $S$
$$(\add S'{\mathcal Z})/S'T'{\mathcal Y}\stackrel S\lra (\add SS'{\mathcal Z})/SS'T'{\mathcal Y}$$
and $T$
$$(\add SS'{\mathcal Z})/SS'T'{\mathcal Y}\stackrel T\lra (\add TSS'{\mathcal Z})/{\mathcal Y}$$
are equivalences of categories, and
$$(\add S'{\mathcal Z})/S'T'{\mathcal Y}\stackrel{TS}{\lra}(\add TSS'{\mathcal Z})/{\mathcal Y}$$
is isomorphic to the natural projection.
\end{enumerate}
\end{Theorem}

In the next corollary we shall try to make the statement of item~\ref{item4} more intelligible.

\begin{Cor}
Under the hypotheses of Theorem~\ref{Greenforadjoints}
we have a commutative diagram
$$
\xymatrix{
(\add S'{\mathcal Z})/S'T'{\mathcal Y}\ar[r]|-S\ar[d]_{\textup{can}}&(\add SS'{\mathcal Z})/SS'T'{\mathcal Y}\ar[dl]|-T\\
(\add TSS'{\mathcal Z})/{\mathcal Y}
}
$$
where $S$ and $T$ are equivalences of categories.
\end{Cor}

\section{Triangulated Generalisation of Auslander Kleiner's Green Correspondence}

We consider now, instead of additive categories,
three  triangulated categories and triangle functors
$$
\xymatrix{
{\mathcal D}\ar@/^/[r]^{S'}&{\mathcal H}\ar@/^/[l]^{T'}\ar@/^/[r]^{S}&{\mathcal G}\ar@/^/[l]^T
}
$$
such that $(S,T)$ and $(S',T')$ are adjoint pairs.

We then replace the additive quotient by Verdier localisation \cite{Verdier}.
\begin{itemize}
\item
For a triangulated category $\mathcal T$ we say that the subcategory
${\mathcal S}$ is thick in $\mathcal T$ if it is triangulated and stable under taking direct summands.
\item
If $\mathcal S$ is thick in $\mathcal T$ then
there is a triangulated category ${\mathcal T}_{\mathcal S}$ together
with a universal functor ${\mathcal T}\ra {\mathcal T}_{\mathcal S}$
annihilating $\mathcal S$, and 'universal' with respect to this property.
\item
For any thick subcategory $\mathcal S$ of $\mathcal T$ we have a canonical functor
$$L_{\mathcal S}:{\mathcal T}/{\mathcal S}\lra {\mathcal T}_{\mathcal S}$$
\end{itemize}

\begin{Theorem} \cite{Greentriangulated} (Green correspondence for triangulated categories)
\label{Greencorrespondencefortriangulated}
Let ${\mathcal D}$, ${\mathcal H}$, ${\mathcal G}$ be three triangulated categories
and let $S,S',T,T'$ be triangle functors
$$
\xymatrix{
{\mathcal D}\ar@/^/[r]^{S'}&{\mathcal H}\ar@/^/[l]^{T'}\ar@/^/[r]^{S}&{\mathcal G}\ar@/^/[l]^T
}
$$
such that $(S,T)$ and $(S',T')$ are adjoint pairs. Let $\epsilon:id_{\mathcal H}\lra TS$
be the unit of the adjunction $(S,T)$. Assume that there is an
endofunctor $U$ of $\mathcal H$ such that $TS=id_{\mathcal H}\oplus U$, denote by
$p_1:TS\lra id_{\mathcal H}$
the projection, and suppose that $p_1\circ\epsilon$ is an isomorphism.

Let ${\mathcal Y}$ be a thick subcategory of $\mathcal H$,
put ${\mathcal Z}:= (US')^{-1}({\mathcal Y})$, and suppose that
each object of $TSS'T'{\mathcal Y}$ is a direct factor of an object of
${\mathcal Y}$.
\begin{enumerate}
\item\label{(1)}
Then $S$ and $T$ induce triangle functors fitting into the commutative diagram
$$
\xymatrix{
({\thick (S'{\mathcal Z})})_{(\thick (S'T'{\mathcal Y}))}\ar[r]|-{S_Z}\ar[d]_{\text{can}}&
({\thick (SS'{\mathcal Z})})_{(\thick (SS'T'{\mathcal Y}))}\ar[dl]|-{T_Z}\\
(\thick (S'{\mathcal Z}))_{{\mathcal Y}}}
$$
of Verdier localisations.
\item\label{(2)}
There is an additive functor $S_\thick$, induced by $S$, and an additive
functor $T_\thick$ induced by $T$,
making  the diagram
$$\xymatrix{(S'{\mathcal Z})/{ (S'T'{\mathcal Y})}\ar[r]^-{\pi_1}\ar[d]^S&
({\thick (S'{\mathcal Z})})/{\thick(S'T'{\mathcal Y})}\ar[d]^{S_\thick}\\
(SS'{\mathcal Z})/{( SS'T'{\mathcal Y})}\ar[r]^-{\pi_2}\ar[d]^T&(\thick (SS'{\mathcal Z}))/{\thick( SS'T'{\mathcal Y})}\ar[d]^{T_\thick}\\
S'{\mathcal Z}/{\mathcal Y}\ar[r]^-{\pi_3}&\thick(S'{\mathcal Z})/\thick({\mathcal Y})
}$$
commutative. Moreover,
the restriction to the respective images of $\pi_1$, respectively $\pi_2$, respectively $\pi_3$
of functors $S_\thick$ and $T_\thick$ on the right
is an equivalence.
\item\label{(3)}
$S$ and $T$ induce equivalences $S_L$ and $T_L$ of additive categories fitting into the commutative diagram
$$
\xymatrix{({\thick (S'{\mathcal Z})})_{(\thick (S'T'{\mathcal Y}))}\ar[rr]|-S\ar@/_7pc/[ddd]_{\text{can}}&&
({\thick (SS'{\mathcal Z})})_{(\thick (SS'T'{\mathcal Y}))}\ar@/^2pc/[dddll]|-T\\
L_{S'T'{\mathcal Y}}((S'{\mathcal Z})/{(S'T'{\mathcal Y})})\ar[r]|-{S_L}\ar[d]_{\text{can}}\ar@^{(->}[u]&
L_{SS'T'{\mathcal Y}} ((SS'{\mathcal Z})/{(SS'T'{\mathcal Y})})\ar[dl]|-{T_L}\ar@^{(->}[ur]\\
L_{{\mathcal Y}}((S'{\mathcal Z})/{\mathcal Y})\ar@^{(->}[d]\\
(\thick (S'{\mathcal Z}))_{{\mathcal Y}}&&&
}
$$
where the outer triangle consists of triangulated categories and triangle functors,
and the inner triangle are full
additive subcategories.
\item\label{(4)}
If $S$ and $T$ induce equivalences of additive categories
$$
\xymatrix{
({\thick (S'{\mathcal Z})})/{ \thick(S'T'{\mathcal Y})}\ar[r]|-{\widehat S}\ar[d]_{\text{can}}&
(\thick (SS'{\mathcal Z}))/{\thick( SS'T'{\mathcal Y})}\ar[dl]|-{\widehat T}\\
(\thick (S'{\mathcal Z}))/{\thick{\mathcal Y}}},
$$
then $S_Z$ and $T_Z$ are mutually inverse equivalences of triangulated categories.
\end{enumerate}
\end{Theorem}

The proof of Theorem~\ref{Greencorrespondencefortriangulated} is quite involved and we refer to
\cite{Greentriangulated} for the interested reader.

\medskip

We shall study briefly what our concept will say for what is known as a localisation, respectively colocalisation  sequence.
We follow Murfet~\cite{Murfet} for the definitions of the following concepts.

A localisation sequence is given by three categories
$\mathcal A$, $\mathcal B$, $\mathcal C$
and functors
$$
\xymatrix{
{\mathcal A}\ar[r]|-i&{\mathcal B}\ar@/^/[l]|-j\ar[r]|-e&
{\mathcal C}\ar@/^/[l]|-\ell
}
$$
such that in addition $(j,i)$,  $(\ell,e)$,  are pairs of adjoint functors, and
such that the counit $\ell e\lra 1$ is an isomorphism, such that the unit $1\lra ij$ is an isomorphism,
and $i$ is the kernel of $e$.

A  colocalisation sequence is given by three categories
$\mathcal A$, $\mathcal B$, $\mathcal C$
and functors
$$
\xymatrix{
{\mathcal A}\ar[r]|-i&{\mathcal B}\ar@/_/[l]|-k\ar[r]|-e&
{\mathcal C}\ar@/_/[l]|-r
}
$$
such that in addition $(k,i)$,  $(e,r)$,  are pairs of adjoint functors, and
such that the unit $1\lra re$ is an isomorphism, such that the counit $ki\lra 1$ is
an isomorphism, and such that $i$ is the
kernel of $e$.

A recollement diagram is given by three categories $\mathcal A$, $\mathcal B$, $\mathcal C$
and functors
$$
\xymatrix{
{\mathcal A}\ar[r]|-i&{\mathcal B}\ar@/^/[l]|-j\ar@/_/[l]|-k\ar[r]|-e&
{\mathcal C}\ar@/^/[l]|-\ell\ar@/_/[l]|-r
}
$$
such that in addition $(j,i)$, $(i,k)$, $(\ell,e)$, $(e,\rho)$ are pairs of adjoint functors, and
such that the unit $1\lra re$ and the counit $\ell e\lra 1$ are isomorphisms and $i$ is the
kernel of $e$.

Consider now the case when
$$
\xymatrix{
{\mathcal D}\ar@/^/[r]^{S'}&{\mathcal H}\ar@/^/[l]^{T'}\ar@/^/[r]^{S}&{\mathcal G}\ar@/^/[l]^T
}
$$
is a localisation sequence. Further we need to assume that the unit $\epsilon:1_{\mathcal H}\lra TS$
has the property that there is an endofunctor $U$ with $TS=I\oplus U$ and the projection $p_1$
onto the first component composes to the identity with $\epsilon.$ Since $SS'=0$ the right end of the triangle
in
Theorem~\ref{Greencorrespondencefortriangulated} is $0$.
Hence, in this situation the statement of Theorem~\ref{Greencorrespondencefortriangulated}
is void.

\end{document}